\documentclass{article}
\usepackage{amsmath}
\usepackage{amssymb}

\usepackage{theorem}
\numberwithin{equation}{section}

\newtheorem{theorem}{Theorem}[section]

\newtheorem{lemma}[theorem]{Lemma}
\newtheorem{definition}{Definition}[section]

\begin{document}
\title{Global well - posedness for the defocusing, cubic, nonlinear wave equation in three dimensions for radial initial data in $\dot{H}^{s} \times \dot{H}^{s - 1}$, $s > \frac{1}{2}$}

\author{Benjamin Dodson}

\maketitle

\noindent \textbf{Abstract:} In this paper we study the defocusing, cubic nonlinear wave equation in three dimensions with radial initial data. The critical space is $\dot{H}^{1/2} \times \dot{H}^{-1/2}$. We show that if the initial data is radial and lies in $(\dot{H}^{s} \times \dot{H}^{s - 1}) \cap (\dot{H}^{1/2} \times \dot{H}^{-1/2})$ for some $s > \frac{1}{2}$, then the cubic initial value problem is globally well - posed. The proof utilizes the I - method, long time Strichartz estimates, and local energy decay. This method is quite similar to the method used in \cite{D2}.

\section{Introduction}
\noindent In this paper we study the defocusing, cubic wave equation

\begin{equation}\label{wave}
u_{tt} - \Delta u = F(u) = -u^{3}, \hspace{5mm} u(0,x) = u_{0}, \hspace{5mm} u_{t}(0,x) = u_{1}, \hspace{5mm} u : \mathbf{R} \times \mathbf{R}^{3} \rightarrow \mathbf{R}.
\end{equation}

\noindent A solution to $(\ref{wave})$ actually produces a family of solutions due to scaling. Indeed, if $u$ solves $(\ref{wave})$ with initial data $(u(0), u_{t}(0))$ then for any $\lambda > 0$,

\begin{equation}\label{scaling}
u(t,x) \mapsto \lambda u(\lambda t, \lambda x),
\end{equation}

\noindent is a solution to $(\ref{wave})$ with initial data $(\lambda u(0, \lambda x), \lambda^{2} u_{t}(0, \lambda x))$. $(\ref{scaling})$ preserves the $\dot{H}^{1/2}(\mathbf{R}^{3})$ norm of $u$ and the $\dot{H}^{-1/2}(\mathbf{R}^{3})$ norm of $u_{t}$, and thus $(\ref{wave})$ is called $\dot{H}^{1/2}$ - critical.\vspace{5mm}

\noindent Study of dispersive partial differential equations with initial data lying in the critical Sobolev space is currently an important topic of research. \cite{Lindblad} and \cite{LS} proved a sharp counterexample to well - posedness for data lying in a Sobolev space less regular than the critical Sobolev space. See \cite{ChrColTao2} for similar results for a number of dispersive equations, including $(\ref{wave})$.\vspace{5mm}

\noindent On the other hand, positive results have been obtained for a number of initial value problems with initial data lying in the critical Sobolev space. \cite{LS} proved a local well - posedness result for $(\ref{wave})$ with initial data in $\dot{H}^{1/2} \times \dot{H}^{-1/2}$.

\noindent For the energy - critical, defocusing wave equation in three dimensions (cubic), global existence of smooth, radially symmetric solutions was proved in \cite{Struwe}. \cite{Grillakis1} extended this result to the general case. \cite{SS} extended this result to dimensions $3 \leq d \leq 7$. Global well - posedness for initial data lying in the energy space was proved by \cite{SS1} and \cite{Kapitanski}.\vspace{5mm}

\noindent \textbf{Remark:} This question has also been completely worked out for the defocusing energy - critical (quintic) Schr{\"o}dinger equation (\cite{B4}, \cite{CKSTT4}), and the defocusing, mass - critical Schr{\"o}dinger equation (\cite{D1}, \cite{KVZ}). In each case scattering has also been proved.\vspace{5mm}

\noindent \textbf{Remark:} The above discussion was not intended to be a complete discussion of defocusing energy - critical and mass - critical problems. For one thing, discussion of dimensions other than $d = 3$ was omitted entirely. Discussion of the focusing problem, see for example \cite{KM3}, was also completely omitted.\vspace{5mm}

\noindent What unites the energy - critical wave equation, the energy - critical Schr{\"o}dinger equation, and the mass - critical Schr{\"o}dinger equation is the presence of a conserved quantity that controls the critical Sobolev norm. For example, if $u$ solves the wave equation

\begin{equation}
u_{tt} - \Delta u = -|u|^{p} u,
\end{equation}

\noindent then the energy

\begin{equation}
E(u(t)) = \frac{1}{2} \int_{\mathbf{R}^{d}} |\nabla u(t,x)|^{2} dx + \frac{1}{2} \int |\partial_{t} u(t,x)|^{2} dx + \frac{1}{p + 2} \int |u(t,x)|^{p + 2} dx
\end{equation}

\noindent is conserved. Therefore for $(\ref{wave})$ the energy is given by

\begin{equation}
E(u(t)) = \frac{1}{2} \int |\nabla u(t,x)|^{2} dx + \frac{1}{2} \int |u_{t}(t,x)|^{2} dx + \frac{1}{4} \int |u(t,x)|^{4} dx = E(u(0)).
\end{equation}

\noindent However, there is no known conserved quantity that controls $\| u(t) \|_{\dot{H}^{1/2}(\mathbf{R}^{3})}$ or $\| u_{t}(t) \|_{\dot{H}^{-1/2}(\mathbf{R}^{3})}$. In fact this is the only obstacle to proving global well - posedness and scattering for $(\ref{wave})$ with radial data.

\begin{theorem}\label{t1}
Suppose $u$ solves $(\ref{wave})$ on an interval $I$, $I$ is the maximal interval of existence of the solution, and

\begin{equation}
\| u \|_{L_{t}^{\infty} \dot{H}^{1/2}(I \times \mathbf{R}^{3})} + \| u_{t} \|_{L_{t}^{\infty} \dot{H}^{-1/2}(I \times \mathbf{R}^{3})} < \infty.
\end{equation}

\noindent Then $u$ is global, that is $I = \mathbf{R}$, and $u$ scatters to a free solution both forward and backward in time.
\end{theorem}

\noindent \emph{Proof:} See \cite{DL}.\vspace{5mm}

\noindent The definitions of well - posedness and scattering that are used here are the standard definitions.

\begin{definition}[Well - posedness]\label{d1.2}
The initial value problem $(\ref{wave})$ is well - posed on an open interval $I \subset \mathbf{R}$, $0 \in I$, for $(u_{0}, u_{1}) \in (\dot{H}^{s} \cap \dot{H}^{1/2}) \times (\dot{H}^{s - 1} \cap \dot{H}^{-1/2}) = X$ if

\begin{enumerate}
\item $(\ref{wave})$ has a unique solution $u$ lying in $C_{t}^{0}(I ; X)$,

\item The solution satisfies the Duhamel formula

\begin{equation}\label{1.8}
(u(t), u_{t}(t)) = S(t)(u_{0}, u_{1}) - \int_{0}^{t} S(t - \tau)(0, u^{3}) d\tau,
\end{equation}

\noindent where $S(t)(f, g)$ is the solution operator to the linear wave equation $u_{tt} - \Delta u = 0$, $u(0,x) = f(x)$, $u_{t}(0,x) = g(x)$.

\item For any compact $J \subset I$, the map $(u_{0}, u_{1}) \mapsto L_{t,x}^{4}(J \times \mathbf{R}^{3})$ is continuous.
\end{enumerate}

\noindent $(\ref{wave})$ is said to be globally well - posed if $I = \mathbf{R}$.

\end{definition}

\begin{definition}[Scattering]\label{d1.3}
A global solution to $(\ref{wave})$ with initial data $(u_{0}, u_{1}) \in X$ is said to scatter forward in time to some $(u_{0}, u_{1})^{+} \in X$ if

\begin{equation}\label{1.9}
 \lim_{t \rightarrow +\infty} \| (u(t), u_{t}(t)) - S(t)(u_{0}, u_{1})^{+} \|_{X} = 0.
\end{equation}

\noindent Analogously, $u$ is said to scatter backward in time to some $(u_{0}, u_{1})^{-} \in X$ if

\begin{equation}\label{1.9.1}
 \lim_{t \rightarrow -\infty} \| (u(t), u_{t}(t)) - S(t)(u_{0}, u_{1})^{-} \|_{X} = 0.
\end{equation}

\noindent $(\ref{wave})$ is said to be scattering for initial data lying in a certain set if for each $(u_{0}, u_{1})$ lying in that set there exists $(u_{0}, u_{1})^{+}$ and $(u_{0}, u_{1})^{-}$ such that $(\ref{1.9})$ and $(\ref{1.9.1})$ hold, and furthermore, the maps $(u_{0}, u_{1}) \mapsto (u_{0}, u_{1})^{+}$ and $(u_{0}, u_{1}) \mapsto (u_{0}, u_{1})^{-}$ are continuous as functions of $(u_{0}, u_{1})$.
\end{definition}

\noindent For a number of focusing, dispersive partial differential equations, there exist solutions with bounded critical Sobolev norm which fail to be global or fail to scatter. This phenomenon is called type two blow up. Excluding type two blowup, such as in \cite{DL}'s proof of theorem $\ref{t1}$, utilizes concentration compactness arguments.\vspace{5mm}

\noindent These arguments are very similar to arguments used to prove global well - posedness and scattering for energy critical wave and Schr{\"o}dinger equations, and mass - critical Schrodinger equations. In fact, given a conserved quantity that controls the critical Sobolev norm, all that is left is to exclude type two blowup. Thus, when \cite{KM2} proved global well - posedness and scattering for the cubic nonlinear Schr{\"o}dinger equation with bounded $\dot{H}^{1/2}(\mathbf{R}^{3})$, this introduced a number of techniques that were very instrumental in the proofs of energy - critical and mass - critical scattering results.\vspace{5mm}

\noindent To the author's knowledge there are no known methods for proving global well - posedness and scattering for dispersive equations without either assuming the existence of a quantity that conserves the critical Sobolev norm or in fact having such a quantity.\vspace{5mm}

\noindent In this paper we utilize the I - method to prove that for any $s > \frac{1}{2}$ the $\dot{H}^{s} \times \dot{H}^{s - 1}$ norm of $(u(t), u_{t}(t))$ is bounded on any finite compact subset of $\mathbf{R}$. This is enough to prove global well - posedness. 

\begin{theorem}[Main theorem]\label{theorem}
$(\ref{wave})$ is globally well - posed for any radial initial data $(u(0), u_{t}(0)) = (u_{0}, u_{1}) \in \dot{H}^{s}(\mathbf{R}^{3}) \times \dot{H}^{s - 1}(\mathbf{R}^{3}) \cap \dot{H}^{1/2} \times \dot{H}^{-1/2}$, $s > \frac{1}{2}$.
\end{theorem}
 
\noindent The I - method has its roots in the Fourier truncation method. The Fourier truncation method was introduced by \cite{B2} for the cubic nonlinear Schr{\"o}dinger equation and by \cite{KPV} for $(\ref{wave})$, proving $(\ref{wave})$ is globally well - posed for $u(0) \in \dot{H}^{s}(\mathbf{R}^{3}) \cap L^{4}(\mathbf{R}^{3})$, $u_{t}(0) \in \dot{H}^{s}(\mathbf{R}^{3})$, $s > \frac{3}{4}$. See also \cite{BC} and \cite{GP}.\vspace{5mm}

\noindent The I - method is an improvement over the Fourier truncation method. For example \cite{CKSTT1} was able to improve the results of \cite{B2} for the nonlinear Schr{\"o}dinger equation. On the wave equation side, \cite{Roy1} extended the results of \cite{CKSTT1} to the inhomogeneous Sobolev spaces $H_{x}^{s}(\mathbf{R}^{3}) \times H_{x}^{s - 1}(\mathbf{R}^{3})$ for $s > \frac{13}{18}$ in the general case and to $s > \frac{7}{10}$ if $u$ has radial symmetry.\vspace{5mm}

\noindent \textbf{Remark:} Inhomogeneous spaces were considered to avoid technical complications at low frequencies.\vspace{5mm}

\noindent Perhaps more importantly, \cite{KTwm} proved a well - posedness result which was technically unattainable via the Fourier truncation method. See \cite{D2} for a more detailed discussion of the history of the I - method.\vspace{5mm}

\noindent To prove our result we make use of the long - time Strichartz estimates. The long time Strichartz estimates were introduced in \cite{D2} and were actually inspired in large part by the linear - nonlinear decomposition of \cite{Roy1}. Basically, the idea is that if $u$ solves $(\ref{wave})$ on an interval $[0, T]$, on which we have some a priori bound on the $\| u(t) \|_{L_{t}^{\infty} \dot{H}^{s}([0, T] \times \mathbf{R}^{3})}$ norm for some $s > \frac{1}{2}$, then we can show that at high frequencies, the solution $u$ is dominated by the free evolution from initial data $(u(0), u_{t}(0))$.\vspace{5mm}

\noindent We then take the usual modified energy

\begin{equation}
E(Iu(t)) = \frac{1}{2} \int |\nabla Iu(t,x)|^{2} dx + \frac{1}{2} \int |Iu_{t}(t,x)|^{2} dx + \frac{1}{4} \int |Iu(t,x)|^{4} dx,
\end{equation}

\noindent where $I$ is a smoothing Fourier multiplier

\begin{equation}
I : \dot{H}^{s}(\mathbf{R}^{3}) \rightarrow \dot{H}^{1}(\mathbf{R}^{3}), \hspace{5mm} I : \dot{H}^{s - 1}(\mathbf{R}^{3}) \rightarrow L^{2}(\mathbf{R}^{3}).
\end{equation}

\noindent Direct computation shows that $\frac{d}{dt} E(Iu(t))$ is a quadrilinear integral operator on $u$ that has at least two terms at high frequencies. Using the long time Strichartz estimates, we can then show that the integral of $\frac{d}{dt} E(Iu(t))$ over the interval $[0, T]$ is small, which in turn implies that $E(Iu(t))$ is pretty close to $E(Iu(0))$. Meanwhile, an a priori upper bound on $E(Iu(t))$ gives us good control over $\| u(t) \|_{\dot{H}^{s}}$, allowing us to make a bootstrap argument which proves theorem $\ref{theorem}$.\vspace{5mm}

\noindent This argument is extremely similar to the scattering argument in \cite{D2}. There are two main reasons we do not prove scattering here. The first is the lack of an interaction Morawetz estimate for the wave equation, unlike the interaction Morawetz estimate for the nonlinear Schr{\"o}dinger equation in \cite{CKSTT2}. The second is that the $L^{2}$ norm of $u$ is not conserved for the nonlinear wave equation $(\ref{wave})$, as it is for the nonlinear Schr{\"o}dinger equation. Observe that \cite{D2} assumed that the initial data lay in $L^{2}(\mathbf{R}^{d})$.

\section{Linear estimates for the wave equation}
In this section we prove some Strichartz - type estimates on solutions to linear wave equations that will be needed in the proof of theorem $\ref{theorem}$. We begin with a discussion of the Littlewood - Paley partition of unity.

\begin{definition}[Littlewood - Paley partition of unity]
Suppose $\psi \in C_{0}^{\infty}(\mathbf{R}^{3})$ is a radial, decreasing function supported on $|x| \leq 2$, $\psi = 1$ on $|x| \leq 1$. Then for any $N$ we define the Littlewood - Paley projection

\begin{equation}
(P_{N} f)(x) = \mathcal F^{-1} ((\psi(\frac{\xi}{N}) - \psi(\frac{2 \xi}{N})) \hat{f}(\xi))(x),
\end{equation}

\noindent where

\begin{equation}
\mathcal F^{-1}(\hat{f}(\xi))(x) = (2 \pi)^{-3/2} \int e^{ix \cdot \xi} \hat{f}(\xi) d\xi.
\end{equation}

\noindent Also define the operators

\begin{equation}
(P_{\leq N} f)(x) = \mathcal F^{-1} (\psi(\frac{\xi}{N}) \hat{f}(\xi))(x),
\end{equation}

\noindent and $P_{> N} = 1 - P_{\leq N}$.
\end{definition}

\noindent \textbf{Remark:} Since $\psi$ is a $C_{0}^{\infty}(\mathbf{R}^{3})$ function, $P_{N} f$ is the convolution of $f$ with a Schwartz function that is $\lesssim_{l} N^{3} (1 + N |x|)^{-l}$ for any $l \in \mathbf{Z}$.\vspace{5mm}

\noindent Next recall the Strichartz estimates of \cite{Stri}.

\begin{theorem}[Strichartz estimate]\label{tStrichartz}
If $u$ solves $u_{tt} - \Delta u = F$ on an interval $I$, with $t_{0} \in I$, then

\begin{equation}
\| u(t) \|_{L_{t,x}^{4}(I \times \mathbf{R}^{3})} \lesssim \| u(t_{0}) \|_{\dot{H}^{1/2}(\mathbf{R}^{3})} + \| u_{t}(t_{0}) \|_{\dot{H}^{-1/2}(\mathbf{R}^{3})} + \| F \|_{L_{t,x}^{4/3}(I \times \mathbf{R}^{3})}.
\end{equation}
\end{theorem}

\noindent \cite{GV} extended Strichartz estimates to all admissible pairs when $d = 3$. Combining Strichartz estimates with local energy decay:

\begin{theorem}[Strichartz estimates]\label{t2.1}
If $u$ solves the wave equation

\begin{equation}\label{2.1}
u_{tt} - \Delta u = F_{1} + F_{2} + F_{3}, \hspace{5mm} u(0) = u_{0}, \hspace{5mm} u_{t}(0) = u_{1},
\end{equation}

\noindent then

\begin{equation}\label{2.2}
\aligned
\| |\nabla|^{1/2} u \|_{L_{t, x}^{4}(I \times \mathbf{R}^{3})} + (\sup_{R} \frac{1}{R^{1/2}} \| \nabla u \|_{L_{t,x}^{2}(I \times \{ |x| \leq R \})}) \\ + \| |\nabla|^{-1/2} u_{t} \|_{L_{t,x}^{4}(I \times \mathbf{R}^{3})} + (\sup_{R} \frac{1}{R^{1/2}} \| u_{t} \|_{L_{t,x}^{2}(I \times \{ |x| \leq R \})}) \\ \lesssim \| u_{0} \|_{\dot{H}^{1}(\mathbf{R}^{3})} + \| u_{1} \|_{L^{2}(\mathbf{R}^{3})} + \| F_{3} \|_{L_{t}^{1} L_{x}^{2}(I \times \mathbf{R}^{3})} \\ + \| |\nabla|^{1/2} F_{1} \|_{L_{t, x}^{4/3}(I \times \mathbf{R}^{3})} + \sum_{j = -\infty}^{\infty} 2^{j/2} \| F_{2} \|_{L_{t,x}^{2}(I \times \{ 2^{j} \leq |x| \leq 2^{j + 1} \})}.
\endaligned
\end{equation}
\end{theorem}

\noindent \emph{Proof:} Again let $S(t)(u_{0}, u_{1})$ be the solution operator to $(\ref{2.1})$ with $F_{1} = F_{2} = F_{3} = 0$,

\begin{equation}\label{2.3}
\cos(t \sqrt{-\Delta}) u_{0} + \frac{\sin(t \sqrt{-\Delta})}{\sqrt{-\Delta}} u_{1}.
\end{equation}

\noindent By Strichartz estimates and the sharp Huygens principle,

\begin{equation}\label{2.4}
\aligned
\| |\nabla|^{1/2} S(t)(u_{0}, u_{1}) \|_{L_{t, x}^{4}(\mathbf{R} \times \mathbf{R}^{3})} + (\sup_{R} \frac{1}{R^{1/2}} \| \nabla S(t)(u_{0}, u_{1}) \|_{L_{t,x}^{2}(\mathbf{R} \times \{ |x| \leq R \})}) \\ + \| |\nabla|^{-1/2} \partial_{t} S(t)(u_{0}, u_{1}) \|_{L_{t, x}^{4}(\mathbf{R} \times \mathbf{R}^{3})} + (\sup_{R} \frac{1}{R^{1/2}} \| \partial_{t} S(t)(u_{0}, u_{1}) \|_{L_{t,x}^{2}(\mathbf{R} \times \{ |x| \leq R \})}) \\ \lesssim \| u_{0} \|_{\dot{H}^{1}(\mathbf{R}^{3})} + \| u_{1} \|_{L^{2}(\mathbf{R}^{3})}.
\endaligned
\end{equation}

\noindent \textbf{Remark:} For example, if $u_{1}$ is supported on $T \leq |x| \leq T + R$, then $S(t)(0, u_{1}) = 0$ for $|x| \leq R$ when $t \leq T - R$ or $t > T + 2R$. Combining this with conservation of energy implies $(\ref{2.4})$ for $u_{0} = 0$. Meanwhile, $\nabla S(t)(u_{0}, 0) = S(t)(\nabla u_{0}, 0)$. Finally, $\partial_{t} S(t)(u_{0}, 0) = \nabla \cdot S(t)(0, \nabla u_{0})$, which completes the proof of $(\ref{2.4})$.\vspace{5mm}

\noindent Then by duality, $(\ref{2.4})$, and the Strichartz estimates of \cite{Stri},

\begin{equation}\label{2.5}
\aligned
\| \nabla \int \frac{\sin(-\tau \sqrt{-\Delta})}{\sqrt{-\Delta}} F(\tau) d\tau \|_{L_{x}^{2}(\mathbf{R}^{3})} + \| \partial_{t} \int \frac{\sin((t - \tau) \sqrt{-\Delta})}{\sqrt{-\Delta}} F(\tau) d\tau \|_{L_{x}^{2}(\mathbf{R}^{3})} \\ \lesssim \| |\nabla|^{1/2} F_{1} \|_{L_{t,x}^{4/3}(I \times \mathbf{R}^{3})} + \sum_{j = -\infty}^{\infty} 2^{j/2} \| F_{2} \|_{L_{t,x}^{2}(I \times \{ 2^{j} \leq |x| \leq 2^{j + 1} \})} + \| F_{3} \|_{L_{t}^{1} L_{x}^{2}(I \times \mathbf{R}^{3})}.
\endaligned
\end{equation}

\noindent Therefore, by the Christ - Kiselev lemma of \cite{ChrKis}, when $u_{0} = u_{1} = 0$,

\begin{equation}\label{2.6}
\aligned
\| |\nabla|^{1/2} u \|_{L_{t, x}^{4}(I \times \mathbf{R}^{3})} + \| |\nabla|^{-1/2} u_{t} \|_{L_{t,x}^{4}(I \times \mathbf{R}^{3})} \lesssim \| u_{0} \|_{\dot{H}^{1}(\mathbf{R}^{3})} + \| u_{1} \|_{L^{2}(\mathbf{R}^{3})} \\ + \| |\nabla|^{1/2} F_{1} \|_{L_{t,x}^{4/3}(I \times \mathbf{R}^{3})} + \sum_{j = -\infty}^{\infty} 2^{j/2} \| F_{2} \|_{L_{t,x}^{2}(I \times \{ 2^{j} \leq |x| \leq 2^{j + 1} \})} + \| F_{3} \|_{L_{t}^{1} L_{x}^{2}(I \times \mathbf{R}^{3})},
\endaligned
\end{equation}

\noindent and

\begin{equation}\label{2.7}
\aligned
 (\sup_{R} \frac{1}{R^{1/2}} \| \int_{0}^{t} \frac{\sin((t - \tau) \sqrt{-\Delta})}{\sqrt{-\Delta}} (F_{1} + F_{3})(\tau) d\tau \|_{L_{t,x}^{2}(\mathbf{R} \times \{ |x| \leq R \})}) \\ + (\sup_{R} \frac{1}{R^{1/2}} \| \partial_{t} \int_{0}^{t} \frac{\sin((t - \tau) \sqrt{-\Delta})}{\sqrt{-\Delta}} (F_{1} + F_{3})(\tau) d\tau \|_{L_{t,x}^{2}(\mathbf{R} \times \{ |x| \leq R \})}) \\
\lesssim \| |\nabla|^{1/2} F_{1} \|_{L_{t,x}^{4/3}(I \times \mathbf{R}^{3})} + \| F_{3} \|_{L_{t}^{1} L_{x}^{2}(I \times \mathbf{R}^{3})}.
\endaligned
\end{equation}

\noindent Therefore it only remains to show

\begin{equation}\label{2.8}
\aligned
(\sup_{R > 0} \frac{1}{R^{1/2}} \| \nabla \int_{0}^{t} \frac{\sin((t - \tau) \sqrt{-\Delta})}{\sqrt{-\Delta}} F_{2}(\tau) d\tau \|_{L_{t,x}^{2}(I \times \{ x : |x| \leq R \})}) \\
+ (\sup_{R > 0} \frac{1}{R^{1/2}} \| \partial_{t} \int_{0}^{t} \frac{\sin((t - \tau) \sqrt{-\Delta})}{\sqrt{-\Delta}} F_{2}(\tau) d\tau \|_{L_{t,x}^{2}(I \times \{ x : |x| \leq R \})}) \\
 \lesssim \sum_{j = -\infty}^{\infty} 2^{j/2} \| F_{2} \|_{L_{t,x}^{2}(I \times \{ 2^{j} \leq |x| \leq 2^{j + 1} \})}.
\endaligned
\end{equation}

\noindent Finally, if $u_{0} = u_{1} = F_{1} = F_{3} = 0$ and $F_{2}$ is supported on $\{ x : |x| \leq R \}$, then the sharp Huygens principle implies that the supports of

\begin{equation}\label{2.9}
\int_{\tau \in [0, t] \cap [kR, (k + 1)R]} S(t - \tau)(0, F_{2}) d\tau
\end{equation}

\noindent are finitely overlapping. Since H{\"o}lder's inequality implies

\begin{equation}\label{2.10}
\| F_{2} \|_{L_{t}^{1} L_{x}^{2}([k R, (k + 1) R] \times \mathbf{R}^{3})} \lesssim R^{1/2} \| F_{2} \|_{L_{t}^{1} L_{x}^{2}([k R, (k + 1) R] \times \mathbf{R}^{3})},
\end{equation}

%\noindent By H{\"o}lder's inequality, if $|x| \leq R_{1}$ and $F_{2}$ is supported on $r \leq R_{2}$,

\begin{equation}\label{2.11}
\aligned
(\sup_{R > 0} \frac{1}{R^{1/2}} \| \nabla u \|_{L_{t,x}^{2}(I \times \{ x : |x| \leq R \})}) + (\sup_{R > 0} \frac{1}{R^{1/2}} \| u_{t} \|_{L_{t,x}^{2}(I \times \{ x : |x| \leq R \})}) \\ \lesssim \sum_{j = -\infty}^{\infty} 2^{j/2} \| F_{2} \|_{L_{t,x}^{2}(I \times \{ 2^{j} \leq |x| \leq 2^{j + 1} \})}
\endaligned
\end{equation}

\noindent follows from $(\ref{2.7})$. This completes the proof of theorem $\ref{t2.1}$. $\Box$\vspace{5mm}

\noindent \textbf{Remark:} The same argument also implies that if $P_{N}$ is a Littlewood - Paley multiplier,

\begin{equation}\label{Bernstein}
\aligned
(\sup_{R > 0} N \| P_{N} u \|_{L_{t,x}^{2}(\mathbf{R} \times \{ x : |x| \leq R \})}) \lesssim \| u_{0} \|_{\dot{H}^{1}(\mathbf{R}^{3})} + \| u_{1} \|_{L^{2}(\mathbf{R}^{3})} \\ + \| \nabla F_{1} \|_{L_{t}^{2} L_{x}^{1}(I \times \mathbf{R}^{3})} + \sum_{j = -\infty}^{\infty} 2^{j/2} \| F_{2} \|_{L_{t,x}^{2}(I \times \{ 2^{j} \leq |x| \leq 2^{j + 1} \})},
\endaligned
\end{equation}

\noindent with constant independent of $N$.\vspace{5mm}

\noindent We will also utilize the endpoint Strichartz estimate of \cite{KlMa}.

\begin{theorem}[Endpoint Strichartz estimates]\label{tendpoint}
For $u_{0}$, $u_{1}$ radial,

\begin{equation}
\| S_{L}(t)(u_{0}, u_{1}) \|_{L_{t}^{2} L_{x}^{\infty}(\mathbf{R} \times \mathbf{R}^{3})} \lesssim \| u_{0} \|_{\dot{H}^{1}(\mathbf{R}^{3})} + \| u_{1} \|_{L_{x}^{2}(\mathbf{R}^{3})}.
\end{equation}

\noindent Also, by duality, if $F$ is radial,

\begin{equation}
\| \int_{\mathbf{R}} S_{L}(-t) (0, F)(t) dt \|_{L_{x}^{2}(\mathbf{R}^{3})} \lesssim \| F \|_{L_{t}^{2} L_{x}^{1}(\mathbf{R} \times \mathbf{R}^{3})}.
\end{equation}
\end{theorem}

\noindent \emph{Proof:} See \cite{KlMa}. $\Box$

\section{Proof of the main theorem}
\noindent We follow the work of \cite{CKSTT1} and later \cite{Roy1} and \cite{Roy2}, and define the I - operator $I : H^{s} \rightarrow H^{1}$, where $I$ is given by the Fourier multiplier

\begin{equation}
m(\xi) = 
\left\{
	\begin{array}{ll}
		1  & \mbox{if } |\xi| \leq N \\
		\frac{N^{1 - s}}{|\xi|^{1 - s}} & \mbox{if } |\xi| > 2N.
	\end{array}
\right.
\end{equation}

\noindent By the Sobolev embedding theorem,

\begin{equation}\label{Sobolev}
\aligned
E(Iu(0)) \lesssim \| \nabla Iu(0,x) \|_{L^{2}}^{2} + \| Iu_{t}(0,x) \|_{L^{2}}^{2} + \| Iu(0,x) \|_{L^{6}(\mathbf{R}^{3})}^{2} \| u(0,x) \|_{L^{3}(\mathbf{R}^{3})}^{2} \\ \lesssim_{\| u_{0} \|_{\dot{H}^{1/2}}} \| \nabla I u(t,x) \|_{L^{2}}^{2} + \| Iu_{t}(t,x) \|_{L^{2}}^{2}.
\endaligned
\end{equation}

\noindent Therefore,

\begin{equation}
E(Iu(0)) \leq C(\| u_{0} \|_{\dot{H}^{1/2}} + \| u_{1} \|_{\dot{H}^{-1/2}}, \| u_{0} \|_{\dot{H}^{s}} + \| u_{1} \|_{\dot{H}^{s - 1}}) N^{2(1 - s)}.
\end{equation}

\noindent To prove global well - posedness it suffices to prove that for any compact interval $[0, T_{0}] \subset \mathbf{R}$, there exists an $N(T_{0})$ sufficiently large so that

\begin{equation}
E(Iu(t)) \leq 2 C N^{2(1 - s)}.
\end{equation}

%\noindent Indeed, by Duhamel's principle,

%\begin{equation}
%Iu(t) = S_{L}(t) (Iu_{0}, Iu_{1}) + \int_{0}^{t} S_{L}(t - \tau) (0, I F(u))(\tau) d\tau,
%\end{equation}

%\noindent so by Strichartz estimates, the Sobolev embedding theorem, H{\"o}lder's inequality in time, and Bernstein's inequality,

%\begin{equation}
%\aligned
%\| \int_{0}^{t} S(t - \tau) (0, I F(u)(\tau)) d\tau \|_{\dot{H}^{1/2}(\mathbf{R}^{3})} \lesssim |t|^{3/4} \| Iu \|_{L_{t}^{\infty} L_{x}^{4}(I \times \mathbf{R}^{3})}^{3} \\ + |t| \| |\nabla|^{1/2} Iu \|_{L_{t}^{\infty} L_{x}^{3}(I \times \mathbf{R}^{3})} \| (1 - I) u \|%_{L_{t}^{\infty} \dot{H}_{x}^{1/2}(I \times \mathbf{R}^{3})}^{2} \lesssim T_{0} N^{\frac{3}{2}(1 - s)}.
%\endaligned
%\end{equation}

\noindent We prove this with a standard bootstrap argument. Suppose that for some interval $[0, T] \subset [0, T_{0}]$,

\begin{equation}
\sup_{t \in [0, T]} E(Iu(t)) \leq 2 C N^{2(1 - s)}.
\end{equation}

\noindent Then we show that for $N(T_{0})$ sufficiently large,

\begin{equation}
E(Iu(t)) \leq \frac{3}{2} C N^{2(1 - s)},
\end{equation}

\noindent which implies $E(Iu(t)) \leq 2 C N^{2(1 - s)}$ on $[0, T_{0}]$.

\begin{definition}
Let $I$ be the Fourier multiplier with a fixed $N$. For $1 \leq M \leq N$, let

\begin{equation}
\aligned
\mathcal S(M) = \| P_{> M} |\nabla|^{1/2} Iu \|_{L_{t,x}^{4}([0, T] \times \mathbf{R}^{3})} + \| P_{> M} |\nabla|^{-1/2} Iu_{t} \|_{L_{t,x}^{4}([0, T] \times \mathbf{R}^{3})} \\ + \sup_{N^{-1} \leq R \leq 4T_{0}} \frac{1}{R^{1/2}} \| P_{> M} \nabla Iu \|_{L_{t,x}^{2}([0, T] \times \{ |x| \leq R \})} \\
+  \sup_{N^{-1} \leq R \leq 4T_{0}} \frac{1}{R^{1/2}} \| P_{> M} Iu_{t} \|_{L_{t,x}^{2}([0, T] \times \{ |x| \leq R \})} \\ + \sup_{N^{-1} \leq R \leq 4T_{0}} \frac{M}{R^{1/2}} \| P_{> M} Iu \|_{L_{t,x}^{2}([0, T] \times \{ |x| \leq R \})}.
\endaligned
\end{equation}
\end{definition}

\begin{theorem}[Long time Strichartz estimate]\label{tlongtime}
Suppose $E(Iu(t)) \leq 2 N^{2(1 - s)}$ on $[0, T]$. Then there exists a small constant $c(s, \| u_{0} \|_{\dot{H}^{1/2}}, \| u_{1} \|_{\dot{H}^{-1/2}}) > 0$ such that if

\begin{equation}
\ln(N) \gtrsim \frac{1 - s}{c(\frac{1}{2} - s)} + \sqrt{\frac{\ln(T_{0})}{c(\frac{1}{2} - s)}},
\end{equation}

\noindent then

\begin{equation}
\mathcal S(\frac{N}{8}) \lesssim N^{1 - s}.
\end{equation}

\end{theorem}

\noindent \emph{Proof:} For a fixed $0 < T' < T$ define

\begin{equation}
\aligned
\mathcal S'(M) = \| P_{> M} |\nabla|^{1/2} Iu \|_{L_{t,x}^{4}([0, T'] \times \mathbf{R}^{3})} + \| P_{> M} |\nabla|^{-1/2} Iu_{t} \|_{L_{t,x}^{4}([0, T'] \times \mathbf{R}^{3})}  \\ + \sup_{N^{-1} \leq R \leq 4T_{0}} \frac{1}{R^{1/2}} \| P_{> M} \nabla Iu \|_{L_{t,x}^{2}([0, T'] \times \{ |x| \leq R \})} \\ +\sup_{N^{-1} \leq R \leq 4T_{0}} \frac{1}{R^{1/2}} \| P_{> M} Iu_{t} \|_{L_{t,x}^{2}([0, T'] \times \{ |x| \leq R \})} \\ + \sup_{N^{-1} \leq R \leq 4 T_{0}} \frac{M}{R^{1/2}} \| P_{> M} u \|_{L_{t,x}^{2}([0, T'] \times \{ x : |x| \leq R \})}.
\endaligned
\end{equation}

\noindent For a large, fixed constant $C_{1}$, let

\begin{equation}
\mathcal T = \{ T' \in [0, T] : \mathcal S'(\frac{N}{8}) \leq C_{1} N^{1 - s}. \}.
\end{equation}

\noindent It is clear from H{\"o}lder's inequality in time and the uniform bound on $E(Iu(t))$ that $\mathcal T$ is nonempty. Also, by the Lebesgue dominated convergence theorem, $\mathcal T$ is a closed set. Therefore, to prove theorem $\ref{tlongtime}$ it suffices to prove that $\mathcal T$ is open in $[0, T]$.\vspace{5mm}

%\noindent First, by the Lebesgue dominated convergence theorem and finite propagation speed, for any $\epsilon > 0$ there exists $R_{0}(\epsilon)$ sufficiently large such that

%\begin{equation}
%\| S_{L}(t) (u_{0}, u_{1}) \|_{L_{t,x}^{4}(|x| > R_{0} + |t|)} \leq \epsilon.
%\end{equation}

%\noindent Therefore, small - data arguments, the endpoint Strichartz estimates of \cite{KlMa} for radial data, and finite propagation speed,

%\begin{equation}
%\| u \|_{L_{t,x}^{4}(|x| > R_{0} + |t|)} \lesssim \epsilon.
%\end{equation}

%\noindent and

%\begin{equation}\label{small}
%\aligned
%\| \nabla Iu(t) \|_{L_{x}^{2}(\{ x: |x| > R_{0} + |t| \})} + \| Iu_{t} \|_{L_{x}^{2}(\{ x : |x| > R_{0} + |t| \})} \\
%\leq \| \nabla Iu_{0} \|_{L_{x}^{2}(\mathbf{R}^{3})} + \| Iu_{t}(0) \|_{L_{x}^{2}(\mathbf{R}^{3})} + \| \nabla Iu \|_{L_{t}^{\infty} L_{x}^{2}([0, T'] \times \mathbf{R}^{3})} \| u \|_{L_{t,x}^{4}(\{ x : |x| > R_{0} + |t| \})}^{2} \leq N^{1 - s} + C \epsilon^{2} N^{1 %- s}.
%\endaligned
%\end{equation}

%\noindent Since we are concerned about long time well - posedness, we will go ahead and assume that $T_{0} > R_{0}$.\vspace{5mm}

\noindent The radial Sobolev embedding theorem implies a bilinear estimate on $[0, T']$ with $T' \in \mathcal T$.

\begin{lemma}[Bilinear estimate]\label{bilinear}
\noindent For $M \leq N$, if $E(Iu(t)) \leq 2C N^{2(1 - s)}$ on $[0, T']$,

\begin{equation}\label{blspace}
\| (P_{> \frac{M}{8}} \nabla Iu) (P_{< N}u) \|_{L_{t,x}^{2} ([0, T'] \times \{ x : |x| \leq 4T_{0} \} } \lesssim (\ln(T_{0}) + \ln(N))^{1/2} S'(\frac{M}{8}) C^{1/2} N^{1 - s},
\end{equation}

\begin{equation}\label{blfreq}
\| (P_{> \frac{M}{8}} u)(P_{< N} u) \|_{L_{t,x}^{2}([0, T'] \times \{ x : |x| \leq 4 T_{0} \})} \lesssim \frac{1}{M}(\ln(T_{0}) + \ln(N))^{1/2} S'(\frac{M}{8}) C^{1/2} N^{1 - s},
\end{equation}

\noindent and

\begin{equation}\label{biltime}
\| (P_{> \frac{M}{8}} Iu_{t}) (P_{< N}u) \|_{L_{t,x}^{2} ([0, T'] \times \{ x : |x| \leq 4T_{0} \} } \lesssim (\ln(T_{0}) + \ln(N))^{1/2} S'(\frac{M}{8}) C^{1/2} N^{1 - s}.
\end{equation}
\end{lemma}

\noindent \textbf{Remark:} Notice that no $I$ is needed in $(\ref{blfreq})$, which is due to the fact $|\xi| m(\xi)$ is increasing in $|\xi|$ and when $|\xi| \sim M$, $|\xi| m(\xi) \sim M$.\vspace{5mm}

\noindent \textbf{Remark:} The proof only uses the fact that $\| P_{< N} u \|_{L_{t}^{\infty} \dot{H}^{1}([0, T'] \times \mathbf{R}^{3})} \lesssim N^{1 - s}$.\vspace{5mm}

\noindent \emph{Proof:} By definition of $\mathcal S'(\frac{M}{8})$,

\begin{equation}
\| P_{> \frac{M}{8}} \nabla Iu \|_{L_{t,x}^{2}([0, T'] \times \{ |x| \leq \frac{1}{N} \})} \lesssim N^{-1/2} \mathcal S'(\frac{M}{8}),
\end{equation}

\noindent so by the Sobolev embedding theorem $\| P_{< N} u \|_{L^{\infty}} \lesssim N^{1/2} \| \nabla P_{< N} u \|_{L^{2}}$,

\begin{equation}
\aligned
\| (P_{> \frac{M}{8}} \nabla Iu) (P_{< N}u) \|_{L_{t,x}^{2} ([0, T'] \times \{ |x| \leq \frac{1}{N} \} } \\ \lesssim \mathcal S'(\frac{M}{8}) \| \nabla Iu \|_{L_{t}^{\infty} L_{x}^{2}([0, T'] \times \mathbf{R}^{3})} \lesssim C N^{1 - s} \mathcal S'(\frac{M}{8}).
\endaligned
\end{equation}

\noindent Now partition $\{ \frac{1}{N} \leq |x| \leq 4T_{0} \}$ into $\lesssim \ln(N) + \ln(T_{0})$ annuli $\{ x : 2^{j} \leq |x| \leq 2^{j + 1} \}$, where $\frac{1}{N} \leq 2^{j} \leq 4T_{0}$. On each annulus, by definition of $\mathcal S'(\frac{M}{8})$,

\begin{equation}
\| P_{> \frac{M}{8}} \nabla Iu \|_{L_{t,x}^{2}([0, T'] \times \{ 2^{j} \leq |x| \leq 2^{j + 1} \})} \lesssim \mathcal S'(\frac{M}{8}) 2^{j/2},
\end{equation}

\noindent while by the radial Sobolev embedding theorem,

\begin{equation}
2^{j/2} \| Iu \|_{L_{t,x}^{\infty}([0, T'] \times \{ 2^{j} \leq |x| \leq 2^{j + 1} \})} \lesssim \| \nabla Iu \|_{L_{t}^{\infty} L_{x}^{2}([0, T'] \times \mathbf{R}^{3})} \lesssim C N^{1 - s}.
\end{equation}

\noindent The arguments to prove $(\ref{blfreq})$ and $(\ref{biltime})$ are identical. $\Box$\vspace{5mm}

\noindent Then by theorem $\ref{t2.1}$, if $\psi \in C_{0}^{\infty}(\mathbf{R}^{3})$, $\psi(x) = 1$ for $|x| \leq 1$, and $\psi(x)$ is supported on $|x| \leq 2$,

\begin{equation}
\aligned
\mathcal S'(M) \lesssim \| \nabla P_{> M} Iu(0) \|_{L_{x}^{2}(\mathbf{R}^{3})} + \| P_{> M} Iu_{t}(0) \|_{L_{x}^{2}(\mathbf{R}^{3})} \\ + \| |\nabla|^{1/2} I((P_{> \frac{M}{8}} u) (P_{> N} u)^{2}) \|_{L_{t,x}^{4/3}([0, T'] \times \mathbf{R}^{3})} \\
+ \sum_{2^{j} \leq 4 T_{0}} 2^{j/2} \| I P_{> M} (\psi(\frac{x}{2 T_{0}}) (P_{> \frac{M}{8}} u)(P_{< N} u)^{2}) \|_{L_{t,x}^{2}([0, T'] \times \{ x : 2^{j} \leq |x| \leq 2^{j + 1} \})} \\
+ \| I P_{> M} ((1 - \psi(\frac{x}{2 T_{0}})) (P_{> \frac{M}{8}} u)(P_{< N} u)^{2}) \|_{L_{t}^{1} L_{x}^{2}([0, T'] \times \mathbf{R}^{3})}.
%+ \| (1 - \psi(\frac{x}{2 T_{0}}) I(((P_{> \frac{M}{8}} u)(P_{< N} u)^{2})) \|_{L_{t}^{1} L_{x}^{2}([0, T'] \times \mathbf{R}^{3})}.
\endaligned
\end{equation}

\noindent Since $|\xi|^{1/2} m(\xi)$ is increasing as $|\xi| \rightarrow \infty$, $|\nabla|^{1/2} I$ obeys a Leibniz type rule. Therefore, by Bernstein's inequality and the definition of $\mathcal S'$,

\begin{equation}
\aligned
\| |\nabla|^{1/2} I((P_{> \frac{M}{8}} u(t)) (P_{> N} u(t))^{2}) \|_{L_{t,x}^{4/3}(\mathbf{R} \times \mathbf{R}^{3})} \\ \lesssim \| |\nabla|^{1/2} I P_{> \frac{M}{8}} u \|_{L_{t,x}^{4}([0, T'] \times \mathbf{R}^{3})} \| P_{> N} u \|_{L_{t,x}^{4}([0, T'] \times \mathbf{R}^{3})}^{2} \\
+ \| P_{> \frac{M}{8}} u \|_{L_{t,x}^{4}([0, T'] \times \mathbf{R}^{3})} \| |\nabla|^{1/2} I P_{> N} u \|_{L_{t,x}^{4}([0, T'] \times \mathbf{R}^{3})} \| P_{> N} u \|_{L_{t,x}^{4}([0, T'] \times \mathbf{R}^{3})} \\ \lesssim \frac{1}{N^{1/2}} \frac{1}{M^{1/2}} \mathcal S'(\frac{M}{8}) \mathcal S'(N)^{2}.
\endaligned
\end{equation}

\noindent Next, by $(\ref{blfreq})$ and the fact that the Littlewood - Paley kernel is rapidly decreasing, $T_{0} >> 1$ and $N >> 1$, and the Sobolev embedding theorems, both radial and standard,

\begin{equation}
\aligned
\sum_{2^{j} \leq 4 T_{0}} 2^{j/2} \| I P_{> M} \psi(\frac{x}{2T_{0}}) ((P_{> \frac{M}{8}} u)(P_{< N} u)^{2}) \|_{L_{t,x}^{2}([0, T'] \times \{ x : 2^{j} \leq |x| \leq 2^{j + 1} \})} \\
\lesssim \mathcal S'(\frac{M}{8}) \frac{C^{1/2} N^{1 - s}}{M} \sum_{2^{j} \leq 4 T_{0}} \| P_{< N} u \|_{L_{t,x}^{\infty}([0, T'] \times \{ x : 2^{j} \leq |x| \leq 2^{j + 1} \})},
\endaligned
\end{equation}

%\noindent so by the radial Sobolev embedding theorem and standard Sobolev embedding theorem,

\begin{equation}
\lesssim \mathcal S'(\frac{M}{8}) \frac{C N^{2(1 - s)}}{M} (\ln(T_{0}) + \ln(N)).
\end{equation}

\noindent Finally, by the radial Sobolev embedding theorem, Holder's inequality in time, and Bernstein's inequality, 

\begin{equation}
\aligned
 \| I P_{> M} ((1 - \psi(\frac{x}{2 T_{0}})) (P_{> \frac{M}{8}} u)(P_{< N} u)^{2}) \|_{L_{t}^{1} L_{x}^{2}([0, T'] \times \{ x : 2^{j} \leq |x| \leq 2^{j + 1} \})} \\
\lesssim \frac{1}{M} \| |x|^{1/2} P_{< N} u \|_{L_{t}^{\infty} L_{x}^{\infty}([0, T'] \times \mathbf{R}^{3})}^{2} \| Iu \|_{L_{t}^{\infty} \dot{H}^{1}([0, T'] \times \mathbf{R}^{3})} \lesssim \frac{C^{3/2} N^{3(1 - s)}}{M}.
\endaligned
\end{equation}

\noindent Therefore,

\begin{equation}
\mathcal S'(M) \lesssim C^{1/2} N^{1 - s} + \frac{C^{3/2} N^{3(1 - s)}}{M} + \mathcal S'(M) \frac{C N^{2(1 - s)}}{M} + \frac{1}{M^{1/2} N^{1/2}} \mathcal S'(M) \mathcal S'(N)^{2},
\end{equation}

\noindent so, by the bootstrap assumption, $\mathcal S'(\frac{N}{8}) \leq C_{1} N^{1 - s}$, if $M \geq N^{\frac{3}{2} - s}$,

\begin{equation}
\mathcal S'(M) \lesssim C^{1/2} N^{1 - s} + C^{3/2} N^{\frac{3}{2} - 2s} + \mathcal S'(M) C N^{\frac{1}{2} - s} + C_{1}^{2} N^{\frac{3}{4} - \frac{3}{2} s} \mathcal S'(M).
\end{equation}

\noindent Therefore, for some $c > 0$ sufficiently small, for $T_{0}$ large and $N$ satisfying

\begin{equation}
\ln(N) \geq \frac{1 - s}{c(s - \frac{1}{2})} + \sqrt{\frac{\ln(T_{0})}{c(s - \frac{1}{2})}},
\end{equation}

\begin{equation}\label{bound}
\mathcal S'(\frac{N}{8}) \lesssim \mathcal S'(N^{\frac{3}{2} - s}) \ln(T_{0}) N^{c \ln(N) (\frac{1}{2} - s)} + C N^{1 - s}.
\end{equation}

\noindent Then theorem $\ref{tlongtime}$ follows from the base case

\begin{lemma}\label{locallemma}
If $u$ solves $(\ref{wave})$ and $E(Iu(t)) \leq 2 C N^{2(1 - s)}$ on $[0, T']$, then

\begin{equation}\label{local}
\mathcal S'(N^{\frac{3}{2} - s}) \lesssim C N^{2(1 - s)} T_{0}^{1/2}.
\end{equation}
\end{lemma}

\noindent Indeed, plugging $(\ref{local})$ into $(\ref{bound})$,

\begin{equation}
\mathcal S'(\frac{N}{8}) << C_{1} N^{1 - s},
\end{equation}

\noindent and therefore $\mathcal T'$ is both open and closed in $[0, T]$. $\Box$\vspace{5mm}

\noindent \emph{Proof of lemma $\ref{locallemma}$:} Since $E(Iu(t)) \leq 2C N^{2(1 - s)}$ for $t \in [0, T] \subset [0, T_{0}]$,

\begin{equation}
\| Iu \|_{L_{t,x}^{4}([0, T] \times \mathbf{R}^{3})}^{4} \lesssim T_{0} 2C N^{2(1 - s)}.
\end{equation}

\noindent Partition $[0, T]$ into $\lesssim \frac{2}{\eta} C T_{0} N^{2(1 - s)}$ subintervals $I_{j}$ such that $|I_{j}| \leq \frac{\eta}{C N^{2(1 - s)}}$, for some small constant $\eta$.  Then on each interval $\| u_{\leq N} \|_{L_{t,x}^{4}(I_{j} \times \mathbf{R}^{3})} \lesssim \eta^{1/4}$.

\noindent Then by theorem $\ref{tStrichartz}$, $E(Iu(t)) \leq 2C N^{2(1 - s)}$, Bernstein's inequality, and the fact that $|\xi|^{1/2} m(\xi)$ is increasing in $|\xi|$,

\begin{equation}
\aligned
\| |\nabla|^{1/2} Iu \|_{L_{t,x}^{4}(I_{j} \times \mathbf{R}^{3})} \lesssim \| Iu \|_{L_{t}^{\infty} \dot{H}^{1}(I_{j} \times \mathbf{R}^{3})} + \| Iu_{t} \|_{L_{t}^{\infty} L_{x}^{2}(I_{j} \times \mathbf{R}^{3})} \\+ \| |\nabla|^{1/2} Iu \|_{L_{t,x}^{4}(I_{j} \times \mathbf{R}^{3})} \| u \|_{L_{t,x}^{4}(I_{j} \times \mathbf{R}^{3})}^{2}  \lesssim N^{1 - s} + \eta^{1/2} \| |\nabla|^{1/2} Iu \|_{L_{t,x}^{4}(I_{j} \times \mathbf{R}^{3})} \\ + \frac{1}{N} \| |\nabla|^{1/2} Iu \|_{L_{t,x}^{4}(I_{j} \times \mathbf{R}^{3})}^{3}.
\endaligned
\end{equation}

\noindent Then since $N$ is large, $\| |\nabla|^{1/2} Iu \|_{L_{t,x}^{4}(I_{j} \times \mathbf{R}^{3})} \lesssim C^{1/2} N^{1 - s}$, and by Bernstein's inequality,

\begin{equation}
\| u \|_{L_{t,x}^{4}(I_{j} \times \mathbf{R}^{3})} \lesssim \| Iu \|_{L_{t,x}^{4}(I_{j} \times \mathbf{R}^{3})} + \| (1 - I)u \|_{L_{t,x}^{4}(I_{j} \times \mathbf{R}^{3})} \lesssim \eta + C^{1/2} N^{\frac{1}{2} - s} \lesssim \eta.
\end{equation}

\noindent Therefore, by theorem $\ref{t2.1}$ and $E(Iu(t)) \leq 2 C N^{1 - s}$ on $I_{j}$,

\begin{equation}
\aligned
\| |\nabla|^{1/2} Iu \|_{L_{t, x}^{4}(I_{j} \times \mathbf{R}^{2})} + (\sup_{R} R^{-1/2} \| \nabla Iu \|_{L_{t,x}^{2}(I_{j} \times \{ |x| \leq R \})}) \lesssim C^{1/2} N^{1 - s}.
\endaligned
\end{equation}

\noindent This proves lemma $\ref{locallemma}$. $\Box$\vspace{5mm}

%\begin{equation}
%\| P_{> N} u \|_{L_{t,x}^{4}(I_{j} \times \mathbf{R}^{3})} \lesssim N^{-1/2} \| |\nabla|^{1/2} Iu \|_{L_{t,x}^{4}(I_{j} \times \mathbf{R}^{3})}.
%\end{equation}

%\noindent Since $N^{1/2 - s} << 1$, the proof of lemma $\ref{locallemma}$, and thus theorem $\ref{tlongtime}$ is complete. $(\ref{biltime})$ follows from the Strichartz estimates in theorem $\ref{t2.1}$. $\Box$\vspace{5mm}

\noindent Theorem $\ref{tlongtime}$ provides a bound on the growth of $E(Iu(t))$.

\begin{lemma}\label{lgrowth}
For any $t \in [0, T_{0}]$, $E(Iu(t)) \leq \frac{3}{2} C N^{1 - s}$.
\end{lemma}

\noindent \emph{Proof:} Again make a bootstrap argument. Let

\begin{equation}
\mathcal T = \{ T \in [0, T_{0}] : E(Iu(t)) \leq \frac{3}{2} C N^{1 - s} \hspace{5mm} \text{for all} \hspace{5mm} t \in [0, T] \}.
\end{equation}

\noindent Because $E(Iu(0)) = C N^{1 - s}$, $\mathcal T$ is clearly nonempty. Also, since $E(Iu(t))$ is a continuous function of time $\mathcal T$ is closed. Therefore, it only remains to show that $\mathcal T$ is open in $[0, T_{0}]$. Then compute

\begin{equation}\label{derivative}
\frac{d}{dt} E(Iu(t)) = \int (Iu_{t})(t,x) (I(u^{3})(t,x) - (Iu)^{3}(t,x)) dx.
\end{equation}

\noindent Splitting $u = u_{h} + u_{l}$, $u_{l} = P_{< \frac{N}{8}} u$, the Fourier support of $u_{l}$ implies that

\begin{equation}
I(u_{l}^{3}) - (Iu_{l})^{3} = 0.
\end{equation}

\noindent Also,

\begin{equation}
I(u_{l}^{2} P_{< \frac{N}{2}} u) - (Iu_{l})^{2} I P_{< \frac{N}{2}} u = 0,
\end{equation}

\noindent which implies that

\begin{equation}
\aligned
(\ref{derivative}) = 3 \int I P_{h} u_{t}(t,x) (I(u_{l}^{2} u_{h})(t,x) - (Iu_{h}(t,x))(Iu_{l}(t,x))^{2}) dx \\ + O(\int Iu_{t}(t,x) (I(u_{h}^{2} u)(t,x) - (Iu_{h})^{2}(t,x) Iu(t,x)) dx).
\endaligned
\end{equation}

\noindent Then by theorem $\ref{tlongtime}$ and lemma $\ref{bilinear}$, for $N$ sufficiently large,

\begin{equation}
\aligned
\int_{0}^{T} \int_{|x| \leq 4 T_{0}} I P_{h} u_{t}(t,x) (I(u_{l}^{2} u_{h})(t,x) - (Iu_{h}(t,x))(Iu_{l}(t,x))^{2}) dx dt \\ \lesssim (\ln(N) + \ln(T_{0})) C_{1}^{2} N^{2(1 - s)} \frac{N^{2(1 - s)}}{N} << N^{2(1 - s)}.
\endaligned
\end{equation}

\noindent Meanwhile, by the radial Sobolev embedding theorem, Bernstein's inequality, the fact that $\frac{1}{N} << T_{0}$, and that the Littlewood - Paley kernel of $I$ is rapidly decreasing outside the ball $|x| \lesssim \frac{1}{N}$,

\begin{equation}
\aligned
\int_{0}^{T} \int_{|x| > 4 T_{0}} I P_{h} u_{t}(t,x) (I(u_{l}^{2} u_{h})(t,x) - (Iu_{h}(t,x))(Iu_{l}(t,x))^{2}) dx dt \\ \lesssim \| |x|^{1/2} u_{l} \|_{L_{t,x}^{\infty}([0, T] \times \mathbf{R}^{3})}^{2} \| \nabla Iu \|_{L_{t}^{\infty} L_{x}^{2}([0, T] \times \mathbf{R}^{3})} \| Iu_{t} \|_{L_{t}^{\infty} L_{x}^{2}([0, T] \times \mathbf{R}^{3})} \\
\lesssim C^{4} N^{4(1 - s)} \frac{1}{N} << N^{2(1 - s)}.
\endaligned
\end{equation}

%\noindent Next we prove the Morawetz estimate,

%\begin{theorem}[Morawetz estimate]
%For $N$ satisfying,

%\begin{equation}
%\ln(N) \gtrsim \frac{1 - s}{c(s - \frac{1}{2})} + \sqrt{\frac{\ln(T_{0})}{c(s - \frac{1}{2})}},
%\end{equation}

%\noindent and $E(Iu(t)) \leq 2 C N^{2(1 - s)}$,

%\begin{equation}
%\int_{0}^{T} \int \frac{1}{|x|} |Iu(t,x)|^{4} dx dt + \int_{0}^{T} \int \frac{1}{|x|} |Iu_{t}(t,x)|^{2} dx dt \lesssim N^{2(1 - s)}.
%\end{equation}
%\end{theorem}

%\noindent \emph{Proof:} Define the Morawetz potential

%\begin{equation}
%M_{1}(t) = \int (Iu_{t}(t,x)) \frac{x}{|x|} \cdot \nabla Iu(t,x) dx.
%\end{equation}

%\noindent Then $\sup_{t \in [0, T]} |M_{1}(t)| \lesssim N^{2(1 - s)}$. Now compute

%\begin{equation}
%\aligned
%\frac{d}{dt} M_{1}(t) = -\int \frac{1}{|x|} |Iu_{t}(t,x)|^{2} dx - \int [\frac{\delta_{jk}}{|x|} - \frac{x_{j} x_{k}}{|x|^{3}}] (\partial_{j} Iu(t,x)) (\partial_{k} Iu(t,x)) dx \\
%+ \int \frac{1}{|x|} |\nabla Iu(t,x)|^{2} dx + \frac{1}{2} \int \frac{1}{|x|} |Iu(t,x)|^{4} dx + \int (F(Iu) - IF(u)) \frac{x}{|x|} \cdot \nabla Iu(t,x) dx.
%\endaligned
%\end{equation} 

%\noindent Next let

%\begin{equation}
%M_{2}(t) = \int (Iu_{t}(t,x)) (Iu(t,x)) \frac{1}{|x|} dx.
%\end{equation}

%\noindent This time Hardy's inequality implies $\sup_{t \in [0, T]} |M_{2}(t)| \lesssim N^{2(1 - s)}$. Also compute

%\begin{equation}
%\frac{d}{dt} M_{2}(t) = 
%\end{equation}

\noindent Next, integrating by parts, again by theorem $\ref{tlongtime}$ and the fact that $|\nabla|^{1/2} I$ satisfies the Leibniz - type rule, and Bernstein's inequality

\begin{equation}
\aligned
\int_{0}^{T} \int (I \partial_{t} u_{h})(t,x) (I(u_{h}(t,x)^{3}) - (Iu_{h}(t,x))^{3}) dx dt \\ = \int_{0}^{T} \int |\nabla|^{-1/2} (I\partial_{t} u_{h})(t,x) |\nabla|^{1/2} (I(u_{h}^{3})(t,x) - (Iu_{h}(t,x))^{3}) dx dt
\endaligned
\end{equation}

\begin{equation}
\aligned
\lesssim \| |\nabla|^{-1/2} I P_{h} u_{t} \|_{L_{t,x}^{4}([0, T] \times \mathbf{R}^{3})} \| |\nabla|^{1/2} Iu_{h} \|_{L_{t,x}^{4}(\mathbf{R} \times \mathbf{R}^{3})} \| u_{h} \|_{L_{t,x}^{4}(\mathbf{R} \times \mathbf{R}^{3})}^{2} \\ \lesssim \frac{C_{1}^{4}}{N} N^{4(1 - s)} << N^{2(1 - s)}.
\endaligned
\end{equation}

\noindent Meanwhile, by lemma $\ref{bilinear}$ and the fact that $\| I P_{< N} u_{t} \|_{\dot{H}^{1}} \lesssim N \| Iu_{t} \|_{L^{2}} \lesssim C N N^{1 - s}$,

\begin{equation}
\aligned
\int_{0} \int_{|x| \leq 4 T_{0}} (I \partial_{t} u_{l})(t,x) (I(u_{h}(t,x)^{3}) - (Iu_{h}(t,x))^{3}) dx dt \\ \lesssim \| u_{h} \|_{L_{t,x}^{4}([0, T] \times \mathbf{R}^{3})}^{2} \| u_{h} (P_{< N} Iu_{t}) \|_{L_{t,x}^{2}([0, T] \times \{ x : |x| \leq 4 T_{0} \})}
\endaligned
\end{equation}

\begin{equation}
\lesssim \frac{C_{1}^{2} N^{2(1 - s)}}{N} C_{1}^{2} N^{2(1 - s)} << N^{2(1 - s)}.
\end{equation}

\noindent Finally, by the radial Sobolev embedding theorem and Bernstein's inequality,

\begin{equation}
\aligned
\int_{0}^{T} \int_{|x| > 4T_{0}} (I \partial_{t} u_{l})(t,x) (I(u_{h}(t,x)^{3}) - (Iu_{h}(t,x))^{3}) dx dt \\ \lesssim \| u_{h} \|_{L_{t,x}^{4}([0, T] \times \mathbf{R}^{3})}^{2} \| |x|^{1/2} I \partial_{t} u_{l} \|_{L_{t,x}^{\infty}([0, T] \times \mathbf{R}^{3})} \| u_{h} \|_{L_{t}^{\infty} L_{x}^{2}([0, T] \times \mathbf{R}^{3})}
\endaligned
\end{equation}

\begin{equation}
\lesssim \frac{C_{1}^{2} N^{2(1 - s)}}{N} C^{2} N^{2(1 - s)} << N^{2(1 - s)}.
\end{equation}

\noindent The term

\begin{equation}
\int_{0}^{T} \int (Iu_{t})(t,x) (I(u_{h}^{2} u_{l})(t,x) - (Iu_{h})^{2}(t,x) Iu_{l}(t,x)) dx dt
\end{equation}

\noindent can be treated as an interpolation of terms with the cubic nonlinearity in the form $u_{l}^{2} u_{h}$ with terms in the cubic nonlinearity of the form $u_{h}^{3}$.\vspace{5mm}

\noindent Therefore, $\int_{0}^{T} |\frac{d}{dt} E(Iu(t))| dt << N^{2(1 - s)}$, which implies that $E(Iu(t)) \leq \frac{3}{2} C N^{2(1 - s)}$, so $[0, T]$ is open in $T_{0}$. Therefore, $E(Iu(t)) \leq C N^{2(1 - s)}$ on $[0, T_{0}]$. $\Box$\vspace{5mm}

\noindent \emph{Proof of theorem $\ref{theorem}$:} Theorem $\ref{theorem}$ follows directly from lemma $\ref{lgrowth}$. By Bernstein's inequality,

\begin{equation}
\aligned
\| u_{> N} \|_{L_{t}^{\infty} \dot{H}^{s}([0, T_{0}] \times \mathbf{R}^{3})} + \| \partial_{t} u_{> N} \|_{L_{t}^{\infty} \dot{H}^{s - 1}([0, T_{0}] \times \mathbf{R}^{3})} \\ \lesssim \frac{1}{N^{1 - s}} \| \nabla Iu \|_{L_{t}^{\infty} L_{x}^{2}([0, T_{0}] \times \mathbf{R}^{3})} + \| Iu_{t} \|_{L_{t}^{\infty} L_{x}^{2}([0, T_{0}] \times \mathbf{R}^{3})} \lesssim 1.
\endaligned
\end{equation}

\noindent Also, 

\begin{equation}
\| u_{< N} \|_{L_{t}^{\infty} \dot{H}^{1}([0, T_{0}] \times \mathbf{R}^{3})} + \| u_{< N} \|_{L_{t}^{\infty} L_{x}^{2}([0, T_{0}] \times \mathbf{R}^{3})} \lesssim N^{1 - s}.
\end{equation}

\noindent Interpolating this bound with the trivial bound

\begin{equation}
\| Iu(t) - Iu(0) \|_{L^{2}(\mathbf{R}^{3})} \leq \int_{0}^{t} \| \partial_{t} Iu(\tau) \|_{L^{2}(\mathbf{R}^{3})} d\tau \lesssim T_{0} N^{1 - s},
\end{equation}

\noindent proves that for $T_{0} > 1$,

\begin{equation}
\| Iu(t) \|_{L_{t}^{\infty} \dot{H}^{s}([0, T_{0}] \times \mathbf{R}^{3})} \lesssim T_{0}^{1 - s} N^{1 - s}
\end{equation}

\noindent and

\begin{equation}
\| Iu(t) \|_{L_{t}^{\infty} \dot{H}^{1/2}([0, T_{0}] \times \mathbf{R}^{3})} \lesssim T_{0}^{1/2} N^{1 - s}.
\end{equation}

\noindent Also for $\frac{3}{p} = \frac{7}{2} - s$, by the Sobolev embedding theorem and definition of $I$,

\begin{equation}
\| Iu_{t}(t) - Iu_{t}(0) \|_{\dot{H}^{-1}(\mathbf{R}^{3})} \leq \int_{0}^{t} \| \Delta Iu(\tau) \|_{\dot{H}^{-1}(\mathbf{R}^{3})} d\tau + \int_{0}^{t} \| Iu^{3}(\tau) \|_{L_{x}^{6/5}(\mathbf{R}^{3})} d\tau
\end{equation}

\begin{equation}
\aligned
\lesssim \int_{0}^{t} N^{1 - s} d\tau + \int_{0}^{t} \| u(\tau) \|_{L_{x}^{3}(\mathbf{R}^{3})}^{2} \| P_{< N} u(\tau) \|_{L_{x}^{6}(\mathbf{R}^{3})} \\ + N^{1 - s} \int_{0}^{t} \| P_{> N} u(\tau) \|_{L_{x}^{3p}(\mathbf{R}^{3})}^{3} d\tau \lesssim T_{0} N^{1 - s} + T_{0}^{3/2} N^{3(1 - s)}.
\endaligned
\end{equation}

\noindent Therefore, by interpolation if $t \in [0, T_{0}]$, $T_{0} > 1$,

\begin{equation}
\| u_{t}(t) \|_{\dot{H}^{s - 1}(\mathbf{R}^{3})} \lesssim T_{0}^{\frac{3}{2}(1 - s)} N^{3(1 - s)},
\end{equation}

\noindent and

\begin{equation}
\| u_{t}(t) \|_{\dot{H}^{-1/2}(\mathbf{R}^{3})} \lesssim T_{0}^{3/4} N^{3/2}.
\end{equation}

\noindent Thus the $\dot{H}^{s} \cap \dot{H}^{1/2} \times \dot{H}^{s - 1} \cap \dot{H}^{-1/2}$ norm is uniformly bounded on any compact subset of $\mathbf{R}$. Global well - posedness then follows from the local result of \cite{LS}. $\Box$

\end{document}